\documentclass[reqno]{amsart}
\usepackage{amssymb,amsthm,euscript,amsmath, mathrsfs, amscd, mathabx}
\usepackage{graphicx}
\usepackage{color}
\usepackage{bm}
\usepackage{hyperref}
\usepackage{baskervald}
\usepackage[cal=boondoxo, scr=boondoxo]{mathalfa}
\usepackage{tikz-cd}

\newcommand{\margh}[1]{}

\def\risom{\overset{\sim}{\rightarrow}}

\input xy
\xyoption{all}
\CompileMatrices

\def\ZZ{{\mathbf Z}}

\def\AA{{\mathbf A}}

\def\PP{{\mathbf P}}

\def\RR{{\mathbf R}}

\def\cG{{\mathcal G}}

\def\cB{{\mathcal B}}

\def\cO{{\mathcal O}}

\def\cV{{\mathcal V}}
\def\cM{{\mathcal M}}
\def\cX{{\mathcal X}}

\def\sD{{\mathscr D}}

\def\sG{{\mathscr G}}

\def\sL{{\mathscr L}}
\def\sM{{\mathscr M}}
\def\sO{{\mathscr O}}

\def\sV{{\mathscr V}}

\def\fg{{\mathfrak g}}

\def\Sch{\operatorname{Sch}}

\def\Hom{\operatorname{Hom}}

\def\Gal{\operatorname{Gal}}

\def\Spec{\operatorname{Spec}}
\def\Pic{\operatorname{Pic}}

\def\EH{\operatorname{EH}}

\def\MR{\operatorname{MR}}

\DeclareMathOperator{\cSh}{\mathcal{S\!h}}

\def\md{\operatorname{md}}
\def\LLSD{\operatorname{LLSD}}

\def\ct{\operatorname{ct}}

\newtheorem{lem}{Lemma}[subsubsection]

\makeatletter
\@addtoreset{lem}{section}
\@addtoreset{lem}{subsection}
\@addtoreset{lem}{subsubsection}
\@addtoreset{lem}{paragraph}
\makeatother

\renewcommand{\thelem}{\ifnum\value{subsubsection}>0{\thesubsubsection.\arabic{lem}}\else{\ifnum\value{subsection}>0{\thesubsection.\arabic{lem}}\else{\thesection.\arabic{lem}}\fi}\fi}

\newtheorem{thm}[lem]{Theorem}
\newtheorem{prop}[lem]{Proposition}
\newtheorem{cor}[lem]{Corollary}

\theoremstyle{definition}
\newtheorem{defn}[lem]{Definition}

\newtheorem{notn}[lem]{Notation}

\theoremstyle{remark}
\newtheorem{rem}[lem]{Remark}
\newtheorem{remark}[lem]{Remark}

\newtheorem{data}[lem]{Data}
\numberwithin{equation}{section}

\begin{document}
\title{Universal limit linear series and descent of moduli spaces}
\author{Max Lieblich}
\author{Brian Osserman}
\begin{abstract}
  We introduce a formalism of descent of moduli spaces, and use it to produce 
  limit linear series moduli spaces for families of curves in which the 
  components of geometric fibers may have nontrivial monodromy. We then 
  construct a universal stack of limit linear series over the
  stack of semistable curves of compact type, and produce new results on
  existence of real curves with few real linear series.
\end{abstract}

\address[Brian Osserman]{Department of Mathematics\\One Shields
Ave.\\ University of California\\Davis, CA 95616}
\address[Max Lieblich]{Department of Mathematics\\University of Washington\\
Padelford Hall\\Seattle, WA 98195}
\thanks{The first named author was partially supported by NSF CAREER
grant DMS-1056129 and NSF Standard Grant DMS-1600813 during this
project.  The second named author was partially supported by
Simons Foundation grant \#279151 during the preparation of this
work.}
\maketitle

\tableofcontents

\section{Introduction}

Given a smooth, projective curve $C$ of genus $g$, recall that a $\fg^r_d$
on $C$ is a linear series of dimension $r$ and degree $d$; in particular,
a $\fg^1_d$ is almost the same as a morphism $C \to \PP^1$ of degree $d$,
up to automorphism of the target. On a general complex curve $C$ of even genus
$g$, it is classical that if we set $d=g/2+1$, then the number of $\fg^1_d$s
on $C$ is given by the Catalan number $\frac{1}{d}\binom{2d-2}{d-1}$.
It is natural to wonder: if we consider instead $C$ a real curve,\footnote{Here
we assume that $C$ is `general' in the sense that it still has
$\frac{1}{d}\binom{2d-2}{d-1}$ complex $\fg^1_d$s.} what are the possible
numbers of real $\fg^1_d$s on $C$?
In \cite{os4}, the second author showed that there exist real curves
$C$ such that all $\frac{1}{d}\binom{2d-2}{d-1}$ of the complex $\fg^1_d$s
are in fact real. More recently, Cools and Coppens \cite{c-c1} showed that
there also exist real curves $C$ such that only
$\binom{d-1}{\left\lceil\frac{d-1}{2}\right\rceil}$ of the $\fg^1_d$s on
$C$ are real. In this paper, we prove the following as a consequence of
general machinery.

\begin{cor}\label{cor:lls-real} Let $g$ be even, and set $d=g/2+1$.
Then
\begin{enumerate}
\item if $d$ is odd, there exist real smooth projective curves of genus
$g$ that carry no real $\fg^1_d$s;
\item if $d$ is even,
there exist real smooth projective curves of genus
$g$ that carry exactly
$$\frac{1}{d-1}\binom{d-1}{d/2}$$
real $\fg^1_d$s.
\end{enumerate}
\end{cor}

In an apparently completely different direction, the geometry
of the moduli spaces $\cM_g$ of curves of genus $g$ has been an active 
subject of
research over the past 30 years, focused on questions such as for which
$g$ the space $\cM_g$ is of general type, and more sharply, what one can say
about the codimension-$1$ subvarieties of $\cM_g$ and their associated
cohomology classes. In this vein, Khosla \cite{kh2} has developed machinery
for constructing and analyzing families of effective divisors on
$\overline{\cM}_g$, but Khosla's theory depends on the existence of
a suitable partial compactification of the universal moduli space of
linear series over $\cM_g$. Specifically, a curve of `compact type' is a 
nodal curve such that removing any node
disconnects the curve; see below for details. Khosla requires an extension 
of the moduli space over a subset of $\overline{\cM}_g$
which is slightly larger than the locus consisting of all curves of
compact type. \footnote{Precisely, Khosla needs a moduli space over 
the space of
`treelike curves,' which are almost of compact type, except that
irreducible components are allowed to have self-nodes. The extension
of our results from curves of compact type to arbitrary treelike curves
is expected to be routine, but involves constructions using compactified
Jacobians which are complementary to what we do in the present paper, and
which we do not pursue.} In the 1980's, Eisenbud and Harris \cite{e-h1}
introduced their theory of limit linear series for curves of compact type, 
which suggests what the fibers should be for the universal space Khosla
requires. However, while Eisenbud and Harris were able to do a great deal 
with their theory, constructions of relative moduli spaces of limit linear
series in families of curves remained
quite limited until recently. In \cite{os20} and \cite{o-m1} the second
author and Murray made substantial progress in this direction, and in
the present paper, we settle the question completely.

Let $\overline{\cM}_g^{\ct}$ denote the open substack of
$\overline{\cM}_g$ parametrizing semistable curves of compact type;
here we may take the base to be $\Spec \ZZ$, or an arbitrary scheme.
We then prove the following.

\begin{thm}\label{thm:universal}
For any triple $g,r,d$, there is a proper relative algebraic space
$$\cG^r_d\to\overline{\cM}_g^{\ct}$$ whose fiber over a point $[C]$
is the moduli scheme $G^r_d(C)$, which is the classical moduli scheme of
linear series in the case that $C$ is smooth, and is a moduli scheme of
limit linear series in the case that $C$ is nodal.
\end{thm}

In particular, we obtain by pullback moduli spaces of (limit) linear
series over arbitrary (flat, proper) families of curves of compact type.

The relationship between Corollary \ref{cor:lls-real} and Theorem
\ref{thm:universal} is as follows. First, Corollary \ref{cor:lls-real}
is proved via degeneration to a curve of compact type, using the
theory of limit linear series. But for both results, it is not enough
to consider families of curves in which the (geometric) components of every 
fiber are defined over the base field. On the contrary, we need to 
consider families where the components have nontrivial monodromy,
in either an arithmetic or geometric sense: in the first case, we need to 
use (geometrically) reducible real curves having components which are
exchanged by complex conjugation, while in the second case, even if
we work over an algebraically closed field, families arise
where components in fibers are exchanged by the monodromy of the
family, and hence do not correspond to components of the total space.

In such cases, despite the foundational advances of \cite{os20} and
\cite{o-m1}, it is not clear how to give a direct definition 
of limit linear series, because we do not have enough line bundles on the
total space to create twists having all the necessary multidegrees. 
In this paper we develop a general formalism of descent of moduli spaces 
and apply it to bypass this difficulty and produce the necessary
moduli spaces of (limit) linear series.
The idea is that any time we have moduli spaces defined only
on certain ``good'' families of varieties, provided the moduli spaces are 
sufficiently functorial with respect to pullback under a morphism
of families, we automatically obtain descent data and thus can
descend them under any \'etale covers. Moreover, when the relevant universal 
family of varieties admits \'etale covers which give ``good'' families, we can 
carry out this descent universally for all families, obtaining canonical
moduli spaces even for families which do not admit ``good'' \'etale covers.
We thus ultimately construct moduli spaces via descent without giving
an intrinsic description of their moduli functors, but we show that the
spaces we construct recover the usual moduli functor under any base change
where the prior constructions apply. This is carried out in general
in \S \ref{sec:descent-lemma}. Next, in \S \ref{sec:lls} we give a 
direct intrinsic definition of a
limit linear series functor as generally as we are able to, we verify
that it satisfies the necessary functoriality condition, and we then
apply the machinery of \S \ref{sec:descent-lemma} to prove 
Theorem \ref{thm:universal}, producing the
desired moduli spaces of limit linear series over arbitrary families of
curves (of compact type). Finally, in \S \ref{sec:real} we apply our
construction to real and $p$-adic curves whose components may not be
defined over the base field, and conclude the proof of 
Corollary \ref{cor:lls-real}. In the process, we develop a more general
tool for studying similar enumerative questions on real (or $p$-adic) 
linear series.  To state the result,
we introduce more precise terminology and notation.

The Eisenbud-Harris definition of limit linear series 
induces a natural scheme structure on an individual
curve $C$ of compact type, which we will
denote by $G^{r,\EH}_d(C)$. There is a natural bijection between $G^r_d(C)$
and $G^{r,\EH}_d(C)$ which is often an isomorphism; see Theorem
\ref{thm:lls-agree} below for details.

Let $K$ be a field and $C_0$ a geometrically reduced, geometrically 
connected, projective curve over
$K$ with (at worst) nodal singularities. 
\begin{itemize}
\item We say that $C_0$ is \emph{totally split} over $K$ if every irreducible
component of $C_0$ is geometrically irreducible, and every node of
$C_0$ is rational over $K$.
\item If $C_0$ is totally split, we say that
$C_0$ is \emph{of compact type} over $K$ if its dual graph is a tree, or
equivalently, if every node of $C_0$ is disconnecting, or
equivalently, if the Picard variety $\Pic^0(C_0)$ of line bundles on
$C_0$ having degree $0$ on every component of $C_0$ is
complete. 
\item Finally, if $C_0$ is not necessarily totally split, let
$L/K$ be an extension such that $C_0$ is totally split over $L$. We
say that $C_0$ is \emph{of compact type} over $K$ if its base extension to 
$C_L$ is of compact type over $L$. (One can check that this is independent 
of the choice of $L$.)
\end{itemize}

We will introduce in Definition \ref{def:smoothing-family} below a notion of 
a ``presmoothing family'' of curves, which is roughly a family of curves
of compact type which admits an \'etale cover over which we know how to
define the functor of limit linear series. We then have the following
result.

\begin{thm}\label{thm:lls-real} Let $K$ be either $\RR$ or a $p$-adic field and
  $C_0$ a genus $g$ curve of compact type over $K$. Let 
$L/K$ be a finite Galois extension such that the extension
$(C_0)_L$ is totally split over $L$. Given $r,d$, let
$$\rho:=g-(r+1)(g+r-d)$$
be the Brill-Noether number, and suppose that $\rho=0$, 
that
the Eisenbud-Harris limit linear series moduli scheme $G^{r,\EH}_d((C_0)_L)$
is finite, and that its set of $\Gal(L/K)$-invariant $L$-rational points
has cardinality $n$ and consists entirely of reduced points.

Given a smooth curve $B$ over $K$, a point $b_0\in B(K)$, and a presmoothing family $\pi:C/B$
such that $C_0=C\times_B b_0$ and the other fibers of $\pi$ are smooth
curves, there exists an open neighborhood $U$ of $b_0$ in the
$K$-analytic topology on $B(K)$
such that for any $b \in U \smallsetminus
\{b_0\}$, the $K$-scheme $G^r_d(C_b)$ of linear series on $C_b$ is
finite, and $G^r_d(C_b)(K)$ has exactly $n$ elements, which are also
reduced points.
\end{thm}

Using Theorem \ref{thm:lls-real}, we can consider a degeneration to
a curve with a rational component glued to $g$ elliptic tails. Because
linear series on $\PP^1$ are parametrized by a Grassmannian, and 
imposing ramification gives a Schubert cycle, the existence of smooth curves 
with a prescribed number of $K$-rational linear series can then be 
reduced to the existence of intersections of Schubert cycles in
a suitable Grassmannian with the corresponding number of $K$-rational
points (see Corollary \ref{cor:lls-schubert} below). Applying this
together with a family of examples produced by Eremenko and Gabrielov
\cite{e-g3} leads to our proof of Corollary \ref{cor:lls-real}.

\subsection*{Acknowledgments} We would like to thank Frank Sottile
for helpful conversations, and especially for bringing the results
of \cite{e-g3} to our attention. We also thank Brendan Creutz, Danny
Krashen, and Bianca Viray for helpful conversations. Finally, we thank
the referee for a heroic reading.

\section{A descent lemma}
\label{sec:descent-lemma}

In this section we prove a simple lemma that will be useful in
constructing the universal limit linear series moduli space. We are
motivated by the following question: suppose we have moduli spaces
associated to a certain collection of ``good'' families of varieties in a 
given class; under what conditions do the moduli spaces extend uniquely to all 
families in the class? Roughly, our answer is that we obtain the desired
extension if the ``good'' families are closed under \'etale base change and 
if the class of varieties admits a `final object up to \'etale covers' which
contains ``good'' families. Since the argument is purely formal, we
work more generally, replacing our class of families of varieties by 
an arbitrary site, and the moduli spaces by sheaves of sets on the slice 
categories of the site.

Precisely, fix a
site $S$. Let $T\subset S$ be a full subcategory. Let $\cSh_S$ denote
the stack of sheaves of sets on $S$. Given a 
subcategory $U\subset S$, let $\cSh_S(U)$ denote the category of
Cartesian functors $U\to\cSh_S$ over $S$.
Concretely, an object of $\cSh_S(U)$ is given by
\begin{enumerate}
\item for each $u\in U$, a sheaf $F_u$ on the
slice category $S_{/u}$;
\item for every arrow $f:u\to v$ between
objects of $U$, an isomorphism of sheaves on $S_{/u}$
      $$\beta_{f}:F_u\risom f^{-1}F_v$$
      such that for every further map $g:v\to
w$ we have
      $$\beta_{g\circ f}=\beta_{f}\circ f^{-1}\beta_{g}.$$
\end{enumerate} Arrows in $\cSh_S(U)$ are given by
compatible maps between the sheaves $F_u$.

Then our main lemma is the following.

\begin{lem}\label{lem:descent-daddy}
  Consider the following conditions.
  \begin{enumerate}
  \item Given an object $t$ of $T$ and a covering $\{s_i\to t\}$ in
    $S$, each $s_i$ is in $T$.
  \item There is a small full subcategory $T_0\subset T$ that is
    closed under all finite products and finite fiber products in $S$ such
    that for any object $s$ 
    of $S$, there is a covering $\{s_i\to s\}$ 
    and maps $s_i\to t_i$ to objects of $T_0$.
  \end{enumerate}
  If $T$ satisfies these conditions then the restriction functor
  $r:\cSh_S(S)\to\cSh_S(T)$ is an equivalence of 
  categories. In particular, any Cartesian functor $G: T\to\cSh_S$
  admits an extension 
  $\widetilde G:S\to\cSh_S$, unique up to unique isomorphism.
\end{lem}

In other words, we can define sheaves on $S$ if we can define them on
a suitable ``sieving category''. We will apply this in Section
\ref{sec:lls} to describe the universal sheaf of limit linear series.
See Remark \ref{rem:descent-lemma-2} for the key examples of situations
in which our conditions (1) and (2) are satisfied.

\begin{proof}
  The proof is aided by Proposition \ref{prop:descent} and Proposition
  \ref{prop:descent-baby} below. We first let $T'$ be the
  category of all objects of $S$ that admit maps to objects of $T_0$
  and let $T''=T\cap T'$.
  The pair $T'\subset S$ satisfies the hypotheses of Proposition
  \ref{prop:descent}, so
  the functor $$\cSh_S(S)\to\cSh_S(T')$$ is an equivalence. The pair
  $T''\subset T'$ satisfies the hypotheses of Proposition
  \ref{prop:descent-baby}, whence the functor
  $$\cSh_S(T')=\cSh_{T'}(T')\to\cSh_{T'}(T'')=\cSh_S(T'')$$ is a chain
  of equivalences.
  Finally, the pair $T''\subset T$ also satisfies the hypotheses of 
  Proposition \ref{prop:descent}, 
  so the natural restriction diagram
  $$\cSh_S(T)\to\cSh_S(T'')$$ is an equivalence. Composing the
  diagrams yields the result.

\end{proof}

\begin{prop}\label{prop:descent}
  Consider the following conditions.
  \begin{enumerate}
  \item Given an object $t$ of $T$ and a covering $\{s_i\to t\}$ in
    $S$, each $s_i$ is in $T$.
  \item Any object $s$ of $S$ admits a covering by objects of $T$.
  \end{enumerate}
  If $T$ satisfies these conditions then the restriction functor
  $r:\cSh_S(S)\to\cSh_S(T)$ is an equivalence of 
  categories. In particular, any Cartesian functor $G: T\to\cSh_S$
  admits an extension 
  $\widetilde G:S\to\cSh_S$, unique up to unique isomorphism.
\end{prop}
\begin{proof}
  Given an object $s\in S$, choose a covering $\{t_i\to s\}$ by
  objects of $T$. By assumption the fiber products $t_{ij}:=t_i\times_s t_j$
  and $t_{ij\ell}:=t_i\times_s t_j\times_s t_\ell$
  exist and lie in $T$. Since $\cSh_S\to S$ is a stack, for any two
  objects $F, F'\in\cSh_S(S)$ the diagram
  $$
  \begin{tikzcd}
    \Hom(F(s),
  F'(s))\ar[r]&\prod\Hom(F(t_i),F(t_j))\ar[r, shift left]\ar[r, shift right]&\prod\Hom(F(t_{ij}),
  F'(t_{ij}))
\end{tikzcd}
$$
is exact. This shows that $r$ is fully faithful.

To show that $r$ is essentially surjective, suppose given
$G\in\cSh_S(T)$. Given $s$ and the covering $\{t_i\to s\}$ as above,
define
$\widetilde G(s)$ to be objects in $\prod G(t_i)$ with descent data
for the covering $\{t_i\to s\}$. Since $\cSh_S$ is a stack, we see
that there is a canonical isomorphism $G\to r\widetilde G$, 
and that $\widetilde G(s)$ is canonically invariant with respect
to the choice of covering $\{t_i\to s\}$. The Cartesian property and
uniqueness of $\widetilde G$ follow from the 
Cartesian property of $G$ by the full-faithfulness already established
and the stack property of $\cSh_S$.
\end{proof}
\begin{prop}\label{prop:descent-baby}
  Consider the following conditions.
  \begin{enumerate}
  \item Given an object $t$ of $T$ and a covering $\{s_i\to t\}$ in
    $S$, each $s_i$ is in $T$.
  \item There is small full subcategory $T_0\subset T$ that is closed under all
    finite products and finite fiber products in $S$ such that every
    object $s$ of $S$ admits a 
    map to an object of $T_0$. 
  \end{enumerate}
  If $T$ satisfies these conditions then the restriction functor
  $r:\cSh_S(S)\to\cSh_S(T)$ is an equivalence of 
  categories. In particular, any Cartesian functor $G: T\to\cSh_S$
  admits an extension 
  $\widetilde G:S\to\cSh_S$, unique up to unique isomorphism.
\end{prop}
\begin{proof}
  Given an object $s\in S$, define a cofiltering category $F_s$ as
  follows. The objects of $F_s$ are arrows $s\to t$ with $t$ in $T_0$. A
  morphism from $\alpha:s\to t$ to $\beta:s\to t'$ is a commutative diagram
  $$
  \begin{tikzcd}
    & s\ar[dl, "\alpha"']\ar[dr, "\beta"] & \\
    t\ar[rr] & & t'.
  \end{tikzcd}
  $$
  The category is cofiltering by the assumption that $T_0$ is closed
  under products and fiber products.

  Suppose given $G\in\cSh_s(T)$.
  The functor sending $\alpha:s\to t$ in $C_t$
  to the sheaf $\alpha^{-1} G(t)$ defines a filtering system of
  isomorphisms of sheaves on the slice category $S_{/s}$. We define
  $\widetilde G(s)$ to be the colimit of this functor. The Cartesian
  property of $\widetilde G$ follows immediately from the fact that
  restriction commutes with colimits. Moreover, if $s$ is in $T$ then,
  since $G$ is Cartesian, we get a canonical isomorphism $\widetilde
  G(s)=G(s)$. Finally, if $\widetilde G\in\cSh_S(S)$ is Cartesian, then we must
  have $\widetilde G(s)$ equal to the colimit of the values over
  $F_s$, since each of those is canonically isomorphic to $\widetilde
  G(s)$ via the pullback maps. This shows that the construction
  $G\mapsto\widetilde G$ is an essential inverse to $r$.
\end{proof}

\begin{remark}\label{rem:descent-lemma-2}
  One simple way that the hypotheses of Lemma \ref{lem:descent-daddy} 
  can be satisfied is if $S$ contains a final object, $T$ is closed 
  under coverings (condition 1), and $T$ contains a covering of the final 
  object of $S$ (condition 2, with $T_0$ the collection of all coverings
  of the final object which lie in $T$). A similar situation arises if $S$ 
  is the big \'etale site of a Deligne-Mumford 
  stack $\cM$; here we can similarly require that $T$ contains an \'etale
  covering of $\cM$. In this case, the absolute products required in $T_0$ 
  are simply fibered products over $\cM$.
\end{remark}

\begin{remark}\label{rem:descent-lemma-1}
  The theory developed here 
  will be relevant to us for the following reason: we will
  start with a family of curves $X\to B$. The relative
  moduli space of limit linear series will only be naturally defined
  over certain base-changed families $X_B'\to B'$, for various $B'\to
  B$. These $B'$ are diverse: they include an \'etale covering of $B$,
  but also all regular schemes mapping into the smooth locus of the
  morphism $X\to B$, and all points mapping to $B$. When we extend the 
  moduli space using the descent
  theory developed here, we want it to retain its value on the $B'\to
  B$ where it can already be defined. Thus, we need to interpolate
  between \'etale coverings and various other $B$-schemes. This is
  what the results of this section accomplish.
\end{remark}

\begin{remark}
Concretely, one can realize Proposition \ref{prop:descent} as follows:
suppose given a sheaf $G_t$ for each $t\in T$ and isomorphisms
$\phi_f:G_t|_u\to G_u$ for each $f:u\to t$ in $T$ that satisfy the
cocycle condition. Then there is a unique sheaf $\widetilde G$ on $S$
whose restriction to each $t\in T$ is canonically isomorphic to $G_t$,
up to unique isomorphism.
\end{remark}

\begin{remark}\label{rem:descent-lemma-3}
  Suppose $S$ is the big \'etale site of $\Spec k$ and $T_0$ is the
  category of iterated fiber products of a single Galois extension
  $k\subset L$. In this case an object of $\cSh_S(T)$ inherits a
  Galois descent datum with respect to the extension $k\subset L$, and
  Lemma \ref{lem:descent-daddy} is solved by Galois
  descent. (Pedantic note: $T$ itself could be much larger, but since
  the sheaves are Cartesian functors on the big site, this does not
  disturb the Galois descent problem.)

  As a consequence, under the equivalence $$r:\cSh_{\Spec k}(\Spec
  k)\to\cSh_{\Spec k}(T)$$ we have the usual isomorphism
  $$G(\Spec k)=r(G)(\Spec L)^{\Gal(L/k)}$$
  identifying the global sections of the descended object with the
  Galois-invariant global sections of the sheaf on $T$. (The Galois
  action is induced by the functoriality.)
\end{remark}

\section{Universal limit linear series by descent}\label{sec:lls}

In this section, we recall the fundamental definitions of
limit linear series functors (largely following the ideas of \S 4 of
\cite{os20}), and verify that they are sufficiently canonical to apply
descent theory to construct the universal moduli space. In fact, because we are
interested in the 
universal setting, we will have to address some new technicalities in
defining limit linear series functors, which are treated by our
definition below of `consistent' smoothing families.

\subsection{Definitions}
\label{sec:definitions}

We begin by defining the families of curves over which we can
define a limit linear series functor, and the more general families
over which we will be able to descend the resulting moduli spaces.

\begin{defn}\label{def:smoothing-family} A morphism of schemes
$\pi: X \rightarrow B$ is a \emph{presmoothing family} if:

\begin{enumerate} 
\item[(I)] $B$ is regular and quasicompact; 
\item[(II)] $\pi$ is flat and proper; 
\item[(III)] The fibers of $\pi$ are curves of compact type; 
\item[(IV)] Any point in the singular locus of $\pi$ which is smoothed in the
generic fiber is regular in the total space of $X$.
\end{enumerate} 

If the following additional conditions are satisfied,
we say that $\pi$ is a \emph{smoothing family}:

\begin{enumerate} 
\item[(V)] $\pi$ admits a section; 
\item[(VI)] every node in every fiber of $\pi$ is
a rational point, and every connected component of the non-smooth locus
maps injectively under $\pi$; 
\item[(VII)] every irreducible component of every fiber of
$\pi$ is geometrically irreducible;
\item[(VIII)] for any connected component $Z$ of the non-smooth locus of 
$\pi$, if $\pi(Z) \neq B$ then $\pi(Z)$ is a principal closed subscheme,
and if $\pi(Z)=B$ and $Y,Y'$ are the closed subschemes of $X$ with
$X=Y \cup Y'$ and $Z=Y \cap Y'$, then 
$\sO_Y(Z)|_Z \otimes \sO_{Y'}(Z)|_Z \cong \sO_Z$.
\end{enumerate}
\end{defn} 

Conditions (VI) and (VII) require in particular that every fiber
is totally split.

It follows from the deformation theory of a nodal curve
together with the regularity of $B$ that the image of a connected
component of the non-smooth locus is always locally principal (see
for instance the deformation theory on p.\ 82 of \cite{d-m}), 
so the
condition that it is principal is always satisfied Zariski locally on
the base. For the last condition, the existence of the stated $Y$ and $Y'$
depends in an essential way on the fibers being of compact type, and 
is proved in Proposition \ref{prop:comps-nodes} below, which also shows
that $Z$ maps isomorphically onto its image in $B$. In particular, both
$\sO_Y(Z)|_Z$ and $\sO_{Y'}(Z)|_Z$ can be trivialized Zariski locally on $B$, 
so we see that both cases of condition (VIII) are satisfied Zariski locally.

\begin{prop}\label{prop:comps-nodes} If $\pi:X \to B$ satisfies
conditions (I)-(VII) of a smoothing family,
then every connected component $Z$ of the non-smooth
locus of $\pi$ is regular and maps isomorphically onto its image in $B$, 
and $\pi^{-1}(\pi(Z))$ may be written as $Y_Z
\cup Y'_Z$, where $Y_Z$ and $Y'_Z$ are closed in $\pi^{-1}(\pi(Z))$,
and $Y_Z \cap Y'_Z = Z$.
\end{prop}

\begin{proof} Our hypotheses imply that the map from $Z$ to $B$ is proper, 
injective and unramified, with trivial residue field extensions. It then 
follows that it is a closed immersion, using the same argument as in
Proposition II.7.3 of \cite{ha1}.
Given this, the regularity is standard; see for instance Proposition
2.1.4 of \cite{os20}. 

Finally, it follows from the regularity that $Z$
(and hence its image) are irreducible, so the existence of $Y_Z$ and $Y'_Z$
in the generic fiber of $\pi(Z)$ follows from condition (VII) of a 
smoothing family. Taking the closures in $\pi^{-1}(\pi(Z))$, we obtain
a decomposition into closed subsets $Y_Z$ and $Y'_Z$ with the correct
intersection in the generic fiber. Then Proposition 15.5.3 of \cite{ega43}
implies that $Y_Z$ and $Y'_Z$ are connected in every fiber, which implies
because the fibers are of compact type that $Y_Z \cap Y'_Z$ is likewise
connected in every fiber. On the
other hand, $Y_Z \cap Y'_Z$ is nonregular in $\pi^{-1}(\pi(Z))$, so
the regularity of $\pi(Z)$ implies that in any fiber, $Y_Z\cap Y'_Z$ must be 
contained among the nodes of that fiber, and it must therefore be equal 
to $Z$.
\end{proof}

We have the following key observations on construction of (pre)smoothing 
families: 

\begin{prop}\label{prop:presmooth-smooth} 
Any quasicompact family of curves $\pi:X \to B$ which is smooth over
$\overline{\cM}_g$ and has fibers of compact type is a presmoothing family. 

Given any presmoothing family $\pi:X \to B$, there exists an \'etale cover 
of $B$ such that the resulting base change of $\pi$ is a smoothing family.
\end{prop}

\begin{proof} For the first assertion, conditions (II) and (III) of a
presmoothing family are immediate, while condition
(I) follows from the smoothness of $\overline{\cM}_g$, and condition (IV)
from the smoothness of the total space of the universal curve over 
$\overline{\cM}_g$ (see Theorem 5.2 of \cite{d-m}).

For the second assertion, it is standard that \'etale base change can
produce sections through smooth points on every component, and this in turn 
can be used to ensure that all components of fibers are geometrically 
irreducible; see the argument for Corollary 4.5.19.3 in 
\textbf{(Err$_{\text{IV}}$, 20)} of \cite{ega44}.
As we have already mentioned, condition (VIII)
is satisfied Zariski locally. Finally, since the non-smooth locus of 
$\pi$ is finite and unramified over $B$, if we fix any $b \in B$, 
Proposition 8(b) of \S 2.3 of \cite{b-l-r} implies that after \'etale base 
change, we will have that each connected component of the 
non-smooth locus has a single (necessarily reduced) point in the fiber $X_b$, 
with trivial residue field extension. 
This implies that after possible further Zariski
localization on the base, the same condition will hold on all fibers,
which then gives that condition (VI) is satisfied as well.
\end{proof}

We will ultimately show that limit linear series moduli spaces
can be descended to any presmoothing family. The smoothing families
are almost the families for which we can give a direct definition of
the limit linear series functor, but we will need to impose one
additional condition, of a more combinatorial nature, which will be
satisfied Zariski locally on the base. Before defining the
condition, we need to introduce a few preliminaries. 

\begin{notn} Given a smoothing family $\pi:X \to B$, let
$Z(\pi)$ denote the set of connected components of the non-smooth
locus of $\pi$, and let $Y(\pi)$ denote the corresponding set of
closed subsets arising as $Y_Z$ or $Y'_Z$ in Proposition
\ref{prop:comps-nodes}.  Let $\varpi:Y(\pi) \to Z(\pi)$ denote the
resulting surjective two-to-one map, and given $Y \in Y(\pi)$, let
$Y^c$ denote the other element of $Y(\pi)$ with $\varpi(Y^c)=\varpi(Y)$.
\end{notn}

Thus, for $Y \in Y(\pi)$, if $Z=\varpi(Y)$, we have that
$\pi^{-1}(\pi(Z))$ is the union along $Z$ of $Y$ with $Y^c$.

The following preliminary definition will be of basic
importance in the definition of a limit linear series.

\begin{defn}\label{defn:multidegree} Let $\pi:X \to B$ be a
smoothing family. Given $d>0$, a \emph{multidegree} on $\pi$ of total
degree $d$ is a map $\md:Y(\pi) \to \ZZ$ such that for each $Y \in
Y(\pi)$, we have $\md(Y)+\md(Y^c)=d$.
\end{defn}

Thus, we can think of a multidegree as specifying how the
total degree $d$ is distributed on each `side' of every node. This
translates into a usual multidegree on every fiber of $\pi$ as
follows.

\begin{prop}\label{prop:multidegree} Let $\pi:X \to B$ be a
smoothing family, and $\md$ a multidegree of total degree $d$. For 
every $b \in B$, there is a unique way to assign integers
$\md_b(Y)$ to each component $Y$ of the fiber $X_b$ such that for
every $Z \in Z(\pi)$ with $b \in \pi(Z)$, and each of the two
$Y_Z \in Y(\pi)$ with $\varpi(Y_Z)=Z$, we have
\begin{equation}\label{eq:mdb-cond} \sum_{Y \subseteq
Y_Z \cap X_b} \md_b(Y) = \md(Y_Z).
\end{equation}
\end{prop}

\begin{proof} The formula \eqref{eq:mdb-cond} implies
immediately that we must have $\sum_{Y \subseteq X_b} \md_b(Y)=d$, by
considering any $Z$ and both choices of $Y_Z$ with $\varpi(Y_Z)=Z$. 
We then see that \eqref{eq:mdb-cond} also determines $\md_b(Y)$ as
\begin{equation}\label{eq:mdb} \md_b(Y)=d-\sum_{Z: Z\cap Y\neq
\emptyset} \md(Y_Z),
\end{equation} where each $Y_Z$ is chosen not to
contain $Y$. Indeed, the sum above is obtained by summing over all nodes on 
$Y$ of the total degrees on the `other side' of each node from $Y$, which 
means that the degree on every component other than $Y$ occurs exactly once 
in the sum. Uniqueness follows immediately, and we also claim that if
$\md_b(Y)$ is given by \eqref{eq:mdb}, it necessarily satisfies
\eqref{eq:mdb-cond}. Indeed, the number of components of $Y_Z$
is equal to the number of nodes lying on $Y_Z$ (including $Z$), and 
if we sum the formula of \eqref{eq:mdb} over all $Y$ in a given $Y_Z$,
we will subtract off $\md(Y_{Z'})$ and $\md(Y_{Z'}^c)$ exactly once each
for all $Z'\neq Z$ meeting $Y_Z$, and we will also subtract off
$\md(Y_Z^c)$ exactly once, so we can rearrange the terms of the sum to
obtain $d-\md(Y_Z^c)=\md(Y_Z)$, together with
$d-\md(Y_{Z'})-\md(Y_{Z'}^c)=0$ for each $Z' \neq Z$ lying on $Y_Z$.
\end{proof}

\begin{defn}\label{defn:concentrated} Let $X_0$ be a totally split
nodal curve.  A multidegree on $X_0$ is \emph{concentrated} on a
component $Y$ if for all $Y' \neq Y$, it assigns degree $0$ to $Y'$.

If $\pi:X \to B$ is a smoothing family, and $d>0$, a
multidegree $\md$ on $\pi$ is \emph{uniformly concentrated} if for all
$b \in B$, the induced $\md_b$ is concentrated on some component $Y$
of $X_b$.
\end{defn}

Our condition on smoothing families is the
following, asserting in essence that we have sufficiently many
uniformly concentrated multidegrees.

\begin{defn}\label{defn:consistent} Let $\pi:X \to B$ be a
smoothing family, and $\md_I=\{\md_i\}_{i \in I}$ a finite collection
of uniformly concentrated multidegrees on $\pi$ of some fixed total degree
$d$. We say that $\md_I$
\emph{sufficient} if for every $b \in B$ and every component $Y \subseteq 
X_b$ there exists $i \in I$ such that the restriction of $\md_i$ to
$X_b$ is concentrated on $Y$.  We say the smoothing family $\pi:X \to
B$ is \emph{consistent} if it admits a sufficient collection of
uniformly concentrated multidegrees.
\end{defn}

\begin{rem} Note that for every $b \in B$ and $Y \subseteq
X_b$, there is always a multidegree on $\pi$ which, when ignoring
other fibers of $\pi$, is concentrated on $Y$. But there does not
appear to be any reason to expect that there is necessarily a
multidegree which is uniformly concentrated and concentrated on
$Y$. In principle, there could be some $b' \in B$ such that if a
multidegree is concentrated on $Y$ in $X_b$, it cannot be concentrated
on any component of $X_{b'}$.  

For instance, consider the case that $B$ is a surface, with
the non-smooth locus of $\pi$ consisting of two connected components
$Z_1$, $Z_2$, whose images in $B$ meet at two points $b_1,b_2$. Thus,
the fibers $X_b$ will have three components when $b=b_1$ or $b_2$. We can
always fix choices of $Y_{Z_1},Y_{Z_1}^c$ and $Y_{Z_2},Y_{Z_2}^c$ so that
in $X_{b_1}$, both $Y_{Z_1}$ and $Y_{Z_2}$ consist of two components, with
$Y_{Z_1} \cap Y_{Z_2}$ then yielding the ``middle'' component of $X_{b_1}$. 
It seems as though it could be possible that we then have $Y_{Z_1}$ and
$Y_{Z_2}$ both consisting of only a single component in $X_{b_2}$, so
that $Y_{Z_1} \cap Y_{Z_2}$ is empty in $X_{b_2}$. If this is the case, we 
see that the only choice of $\md$ which is
concentrated on $Y_{Z_1} \cap Y_{Z_2}$ in $X_{b_1}$ is given by setting
$\md(Y_{Z_1})=d$ and $\md(Y_{Z_2})=d$, but the resulting multidegree in
$X_{b_2}$ will not be concentrated on any component (it will have degree
$d$ on the extremal components and degree $-d$ in the middle).
\end{rem}

\begin{rem}
There is not in general a unique minimal
sufficient collection of uniformly concentrated multidegrees: for
instance, if we have that $B$ is one-dimensional, and there are two
points $b_1,b_2 \in B$ such that $X_b$ is smooth if $b \neq b_1,b_2$
and $X_{b_1}$ and $X_{b_2}$ each have one node, then because the
behavior at $X_{b_1}$ and $X_{b_2}$ is independent, there are four
uniformly concentrated multidegrees, and there are two different ways
of choosing two of these four to get a sufficient collection.

On the other hand, while the collection of all
uniformly concentrated multidegrees is canonical, it turns out not to
be large enough in general to allow for the most transparent treatment
of our construction. Indeed, distinct multidegrees may become the same
under base change, so when considering base change of limit linear series
we will naturally be led to allow the possibility that our collection
$\md_I$ has repeated entries.
\end{rem}

\begin{prop}\label{prop:consistent-local} If $\pi:X \to B$ is
a smoothing family then for any $b \in B$, there is a Zariski
neighborhood $U$ of $b$ on which $\pi$ is consistent.

In addition, if $(\md_i)_{i \in I}$ is a sufficient
collection of uniformly concentrated multidegrees, and $\md$ is any
uniformly concentrated multidegree, then there is a Zariski
neighborhood $V$ of $b$ and $i \in I$ such that $\md=\md_i$ on $V$.
\end{prop}

\begin{proof} Given $b \in B$, using the language of
\cite{os20}, we choose $U$ so that $X_b$ meets every node of $\pi$ over
$U$, and furthermore so that $\pi$ is `almost local' over $U$ -- see
Definition 2.2.2 and Remark 2.2.3 of \cite{os20}. Intuitively, this 
means that every fiber over $U$ is naturally a partial smoothing of the
chosen fiber $X_b$. In this situation,
it is clear that a multidegree is uniformly concentrated over $U$ if
and only if its restriction to $X_b$ is concentrated on some
component, and it follows easily that the smoothing family will be
consistent over $U$.

Similarly, if $(\md_i)_{i \in I}$ is sufficient, then
by definition there is some $i$ such that $\md$ agrees with $\md_i$ on
the fiber $X_b$. But then it is clear that if we set $V$ as above so
that $\pi$ is almost local, then agreement on $X_b$ implies agreement
on all of $V$.
\end{proof}

\begin{notn}\label{notn:twists} Given a smoothing family
$\pi:X \to B$, and $Y \in Y(\pi)$ over $Z \in Z(\pi)$, define $\sO^Y$
as follows: if $\pi(Z) \neq B$, then $\sO^Y=\sO_X(Y)$; if $\pi(Z)=B$,
then $\sO^Y$ is the line bundle on $\sO_X$ obtained by gluing
$\sO_Y(-Z)$ to $\sO_{Y^c}(Z)$ along $Z$.
\end{notn} 

In the above, we use that when $\pi(Z)=B$ we have $Y \cap Y^c \cong B$, so 
that the restriction maps induce a canonical isomorphism $\Pic(X) \risom 
\Pic(Y) \times_{\Pic(Z)} \Pic(Y^c)$. We also use condition (VIII) of a 
smoothing family to ensure that the images in $\Pic(Z)$ agree. In this case 
we also obtain that sections
of $\sO^Y$ correspond to pairs of sections of $\sO_Y(-Z)$ and $\sO_{Y^c}(Z)$
which `agree' on $Z$; agreement depends in principle on an isomorphism
$\sO_Y(-Z)|_Z \cong \sO_{Y^c}(Z)|_Z$, but is for instance canonical in the
cases of sections vanishing along $Z$.

The purpose of this construction is as follows. Given a line bundle $\sL$
on $X$, its multidegree is defined by looking at the degrees of the
restrictions to each $Y\in Y(\pi)$. On fibers, this is equivalent data to a 
multidegree in the usual sense of the degree on each component, as described
by Proposition \ref{prop:multidegree}. Now, given $\sL$ of multidegree $\md$,
for any $Y \in Y(\pi)$ we have that the multidegree of $\sL\otimes \sO^Y$ is
obtained from $\md$ by subtracting $1$ from $Y$ and adding $1$ to
$Y^c$. Consequently, twisting by the different $\sO^Y$, we can change
the multidegree of a line bundle from any multidegree to any
other (of the same total degree). 
Given a line bundle $\sL$ on $X$, and a multidegree $\md$, we
denote by $\sL_{\md}$ the twist of $\sL$ having multidegree $\md$,
which is unique up to isomorphism.

\begin{data} Suppose $\pi:X\to B$ is a consistent smoothing
family.  A \emph{limit linear series datum\/} on $\pi$ consists of the
following:
\begin{enumerate}
\item A sufficient collection $\md_I$ of
uniformly concentrated multidegrees on $\pi$ (always assumed of some fixed 
degree).
\item For each $Y\in Y(\pi)$ a section $s_Y$
of $\sO^Y$ vanishing precisely on $Y$.
\item For each $Z\in Z(\pi)$, an
isomorphism $$\theta_Z:\sO^{Y_Z}\otimes\sO^{Y'_Z}\risom\sO_X,$$ where
$Y_Z$ and $Y'_Z$ are the elements of $Y(\pi)$ lying over $Z$.
\item A distinguished element $i_0 \in I$.
\end{enumerate} We will denote a limit series datum on
$\pi$ by $\sD=(\md_I, \{s_Y\}, \{\theta_Z\}, i_0)$. Given a consistent
smoothing family $\pi$, we will write $\LLSD(X/B)$ for the set of
limit linear series data attached to $\pi$.
\end{data}

The purpose of items (2)-(4) of a limit linear series datum is to allow
us to pin down the line bundles $\sL_{\md}$ precisely, together with
maps between them. Specifically, for any given $\md$, there is a unique
(up to reordering) minimal collection of twists by the $\sO^Y$ in order
to get from multidegree $\md_{i_0}$ to multidegree $\md$; we then define
$\sL_{\md}$ to be obtained from $\sL$ by the corresponding tensor product.
Then, because $s_Y$ induces a map from any line bundle to its twist by
$\sO^Y$, and $\theta_Z \circ s_{Y^c}$ induces a map in the other direction,
for each $\md,\md'$ we have that a limit linear series datum induces a
canonical map $\sL_{\md} \to \sL_{\md'}$.

\begin{prop}\label{prop:LLS-family} If $\pi:X\to B$ is
a consistent smoothing family, then we have
$\LLSD(X/B)\neq\emptyset$.
\end{prop}

\begin{proof} The existence of the $\md_I$ and $i_0$ are by
definition, and likewise the $s_Y$ exist by the construction of
$\sO^Y$. Finally, the $\theta_Z$ exist by construction in the case that
$\pi(Z)=B$, and by the hypothesis that $\pi(Z)$ is principal in $B$ in
the case that $\pi(Z) \neq B$.
\end{proof}

We now state the general definition of limit linear series for
families of curves of compact type.

\begin{defn} Let $f:T \to B$ be a $B$-scheme, and write
$\pi':X\times _B T \to T$. Suppose $\sD:=(\md_I, \{s_Y\},
\{\theta_Z\}, i_0)$ is a limit linear series datum. A \emph{$T$-valued
$\sD$-limit linear series\/} of rank $r$ and degree $d$ on $\pi$
consists of
\begin{enumerate}
\item an invertible sheaf $\sL$ of multidegree
$\md_{i_0}$ on $X\times_B T$, together with
\item for each $i\in I$ a rank $(r+1)$
subbundle $\cV_i \subseteq \pi'_* \sL_{\md_i}$
\end{enumerate} satisfying the following condition:
for any multidegree $\md$ on $\pi$ of total degree $d$, the map
\begin{equation}\label{eq:lls-map} \pi'_* \sL_{\md}
\to \bigoplus_{i \in I} \left(\pi'_* \sL_{\md_i}\right)/\cV_i
\end{equation} induced by the $s_Y$ and $\theta_Z$ has
$(r+1)$st vanishing locus equal to all of $T$.
\end{defn}

In the above, to say that $\cV_i$ is a rank-$(r+1)$ subbundle
of $\pi'_* \sL_{\md_i}$ means that it is locally free of rank $r+1$
and the injection into $\pi'_* \sL_{\md_i}$ is preserved under base
change (where on $\sL_{\md_i}$, base change is applied prior to
pushforward).  The $(r+1)$st vanishing locus of \eqref{eq:lls-map} is
a canonical scheme structure on the closed subset of points on which
the map has kernel of dimension at least $r+1$, defined in terms of
perfect representatives of $R\pi'_* \sL_{\md}$ and the $R\pi'_*
\sL_{\md_i}$; see Appendix B.3 of \cite{os20} for details.   

\begin{defn} Two $T$-valued $\sD$-limit linear series $(\sL,(\sV_i)_i)$
and $(\sL',(\sV'_i)_i)$ are \emph{equivalent} if there exists a line
bundle $\sM$ on $T$ and an isomorphism $\sL \risom \sL' \otimes \pi'^*
\sM$ sending each $\sV_i$ to $\sV'_i \otimes \sM$ under the
identification $\pi'_* \left(\sL' \otimes \pi'^* \sM\right) = \pi'_*
\sL' \otimes \sM$.  
\end{defn}

\begin{defn}\label{defn:d-lls}
Given a consistent smoothing family $\pi:X\to B$
with a limit linear series datum $\sD$, the \emph{functor of
$\sD$-limit linear series\/}, denoted $\sG^r_d(X/B,\sD)$, is the
functor whose value on a $B$-scheme $T\to B$ is the set of equivalence
classes of $T$-valued $\sD$-limit linear series of rank $r$ and degree
$d$.
\end{defn}


\subsection{Functoriality}
\label{sec:functoriality}

In the remainder of this section, we
show that the limit linear series functor is sufficiently canonical to
apply our descent machinery.  The observation of fundamental
importance to us is the following, which we will use to show not only that 
our definition is independent of the choices made, but also 
that it is well-behaved under base change.

\begin{prop}\label{prop:change-of-multdeg} Suppose $\pi:X\to
B$ is a consistent smoothing family with a limit linear series datum
$\sD=(\md_I, \{s_Y\}, \{\theta_Z\}, i_0)$. Suppose $(\sL,(\sV_i)_i)$
is a $\sD$-limit linear series of rank $r$, and let $\md$ be a
uniformly concentrated multidegree, not necessarily contained among
the chosen $\md_i$. Then the map \eqref{eq:lls-map} has empty
$(r+2)$nd vanishing locus, so that the kernel is universally a
subbundle of rank $r+1$.

Moreover, if we let $\cV$ denote the kernel of
\eqref{eq:lls-map}, then for any multidegree $\md'$, the map
$$\pi'_* \sL_{\md'} \to \left(\pi'_* \sL_{\md}\right)/\cV \oplus 
\bigoplus_{i \in I} \left(\pi'_*
\sL_{\md_i}\right)/\cV_i$$ has $(r+1)$st vanishing locus equal to all
of $T$.
\end{prop}

\begin{proof} By definition of a limit linear series,
\eqref{eq:lls-map} has vanishing locus equal to all of $T$, so in
order to prove that the kernel is universally a subbundle of rank
$r+1$, it suffices to show that the $(r+2)$nd vanishing locus is
empty, which is the same as saying that there are no points of $t$ at
which the kernel has dimension at least $r+2$.  Given $t \in T$,
denote by $X'_t$ the corresponding fiber of $\pi'$; then the fiber of
\eqref{eq:lls-map} at $t$ is
$$H^0(X'_t,\sL_{\md}|_{X'_t}) \to 
\bigoplus_{i \in I} H^0(X'_t,
\sL_{\md_i}|_{X'_t})/\cV_i|_t.$$ Now, because $\md$ is uniformly
concentrated, its restriction to $X'_t$ is concentrated on some
component of $X'_t$. On the other hand, because $(\md_i)_i$ is assumed
sufficient, this means that $\md$ is equal to some $\md_i$ after
restriction to $X'_t$, so that $H^0(X'_t,\sL_{\md}|_{X'_t}) =
H^0(X'_t, \sL_{\md_i}|_{X'_t})$. But then the kernel of the above map
is contained in $\sV_i|_t$, which is $(r+1)$-dimensional, proving the
first statement.

Next, the second statement can be checked Zariski
locally, so by Proposition \ref{prop:consistent-local}, we may assume
that in fact $\md=\md_i$ for some $i$, and then necessarily $\cV$ is
identified with $\cV_i$ under the resulting isomorphism $\sL_{\md}
\risom \sL_{\md_i}$. Thus, the kernel of the given map is
(universally) identified with the kernel of \eqref{eq:lls-map} in
multidegree $\md'$, so the $(r+1)$st vanishing loci also agree, as
desired (see Proposition B.3.2 of \cite{os20}).
\end{proof}

\begin{cor}\label{cor:lls-independence} Suppose $\pi:X\to B$
is a consistent smoothing family.
\begin{enumerate}
\item Given two limit linear series data $\sD$
and $\sD'$ on $\pi$, there is a canonical
isomorphism $$\tau_{\sD,\sD'}^{X/B}:\sG^r_d(X/B,\sD)\risom\sG^r_d(X/B,\sD').$$
\item For any base change $B'\to B$ preserving
the smoothing family conditions, there is a canonical isomorphism
      $$\beta_{B'/B}^\sD:\sG^r_d(X_{B'}/B',\sD_{B'})\risom\sG^r_d(X/B,\sD)\times_B
      B'.$$
    \end{enumerate}
      By ``canonical'' we mean that for any
triple $\sD, \sD', \sD''$ of limit linear series data we have that
      $$\tau_{\sD,\sD''}^{X/B} = \tau_{\sD',\sD''}^{X/B} \circ
      \tau_{\sD,\sD'}^{X/B},$$ and similarly
for any triple $B''\to B'\to B$ preserving the smoothing family
conditions, we have that
      $$\beta_{B''/B}^\sD=\beta_{B'/B}^\sD|_{B''}\circ\beta_{B''/B'}^{\sD_{B'}}.$$
\end{cor}
\begin{proof} First, if $(\md_i)_i$ and $i_0$ are fixed, the
choices of the $s_Y$ and $\theta_Z$ are only used to determine a
subfunctor of the functor of all tuples of $(\sL,(\cV_i)_{i \in I})$,
so there is a canonical notion of equality of these subfunctors. We
see from the definitions that both the $s_Y$ and $\theta_Z$ are unique
up to $\sO_T^{\times}$-scalar; for the former, where $Y$ surjects onto
$B$ it is crucial that the (fibers of the) support of $s_Y$ is a connected 
curve, as
otherwise $s_Y$ would only be unique up to independent scaling on each
connected component. This implies
that the rank of \eqref{eq:lls-map} is independent of the choices of
$s_Y$ and $\theta_Z$, as desired. Similarly, because isomorphisms of
line bundles are unique up to $\sO_T^{\times}$-scaling, induced
identifications of subbundles of global sections are unique,
independent of the choice of isomorphism. Because twisting provides a
canonical identification between isomorphism classes of line bundles
of multidegree $\md_{i_0}$ and line bundles of any other given
multidegree (of the same total degree), we conclude that the limit
linear series functors associated to different choices of $i_0$ are
canonically isomorphic.

Next, in order to show independence of the choice of
$(\md_i)_i$, we observe that it is enough to construct canonical
isomorphisms when we add a uniformly concentrated multidegree to an
existing collection: indeed, we can then compare the functors for any
two collections by comparing each to the functor obtained from the
union of the collections.  Accordingly, suppose we have a collection
$(\md_i)_i$, and an additional uniformly concentrated multidegree
$\md'$. We define the obvious forgetful map from the functor of limit
linear series associated to $(\md_i)_i \cup (\md')$ to the one
associated to $(\md_i)_i$, noting that the $(r+1)$st vanishing locus
of \eqref{eq:lls-map} can only increase if we drop a concentrated
multidegree from the target. That this forgetful map is an isomorphism
of functors then follows immediately from Proposition
\ref{prop:change-of-multdeg}, as the only possibility is to have $\cV
\subseteq \pi'_* \sL_{\md'}$ be the kernel of \eqref{eq:lls-map} in
multidegree $\md'$.

It remains to address compatibility with base change. Because the Picard
functor is compatible with base change, the main 
complications arise from changes in $Z(\pi)$, which can occur either
because the image of $B'$ may be disjoint from the image of
some $Z \in Z(\pi)$, or because the preimage in $B'$ of some $Z \in Z(\pi)$
may decompose into two or more connected components. This means in 
particular that limit linear series data may not have canonical pullbacks. 
In general,
there are natural maps $Z(\pi') \to Z(\pi)$ and $Y(\pi') \to Y(\pi)$ which 
induce a pullback
on multidegrees and on uniformly concentrated multidegrees, and a
sufficient collection of uniformly concentrated multidegrees pulls back
to a sufficient collection. However, the pullback map on multidegrees
is in general neither injective nor surjective. If we pull back our 
given collection $(\md_i)_i$ together with $i_0$, we find that as above, we 
are comparing two closed subfunctors of a fixed functor: 
the functor of all tuples of $(\sL,(\cV_i)_{i \in I})$. Note however
that when the map $Z(\pi') \to Z(\pi)$ is not injective, we do not have
that each $\sO^{Y'}$ can be chosen to be a pullback of some $\sO^Y$. In
addition, when $Z(\pi') \to Z(\pi)$ is not surjective, we will have that
different multidegrees for $\pi$ become the same on $\pi'$. Thus, for a given
multidegree $\md$ and its pullback $\md'$ to $\pi'$, it is not necessarily
the case that when we apply our construction to $\pi'$ to obtain the map 
\eqref{eq:lls-map} for $\md'$, the result is exactly equal to
the pullback of the map for $\md$. However, we claim that each summand
of the two maps agree up to $\sO_{B'}^{\times}$-scalar, so that the resulting
$(r+1)$st vanishing loci conditions are the same. Indeed, the map from
multidegree $\md$ to some $\md_i$ is obtained by a sequence of twists
(or inverse twists, using the $\theta_{\varpi(Y)}$) by different $\sO^Y$. If
$\varpi(Y)$ remains nonempty and connected in $\pi'$, then the pullback of
$s_Y$ is a valid choice for $s_{Y'}$, and the claim is clear. Similarly,
if $\varpi(Y)$ becomes empty, then the twist by $\sO^Y$ doesn't change the 
multidegree for $\pi'$, and the pullback of $s_Y$ can be used to
trivialize the pullback of $\sO^Y$. Finally, if $\varpi(Y)$ breaks into
distinct connected components, denote these by $Y'_1,\dots,Y'_m$. Then
twisting by $\sO^Y$ pulls back to a composition of twists by the $Y'_i$,
and the pullback of $s_Y$ will agree up to $\sO_{B'}^{\times}$-scalar with the
product of the $s_{Y'_i}$, so we conclude the claim. It follows that 
the condition on the $(r+1)$st vanishing locus for the given $\pi'$ and 
$\md'$ agrees with the pullback of the same condition for $\pi$ and $\md$.

Finally, although not every $\md'$ for $\pi'$ arises by pullback from
an $\md$ -- because a $Z \in Z(\pi)$ may break into multiple connected
components after base change -- we can check equality of the subfunctors
after Zariski localization. If we restrict to almost local open
subsets of $B'$ (as described in the proof of Proposition 
\ref{prop:consistent-local}), then the maps $Z(\pi') \to Z(\pi)$ and
$Y(\pi')\to Y(\pi)$ become
injective, so the induced map on multidegrees is surjective, and we obtain 
the desired equality of functors.
\end{proof}

\begin{remark} Note that nearly all the hypotheses of a
smoothing family are automatically preserved under base change (to a
suitable base): the only one which is not necessarily preserved is the
regularity hypothesis on $X$, and of course this is still preserved if
the base change is \'etale, or if $B'$ is a point. The only other condition
which is not obviously preserved under base change is condition (VIII).
Although scheme-theoretic
image is not in general preserved under base change, in our situation
Proposition \ref{prop:comps-nodes} says
that the map to $B$ from each connected component of the non-smooth locus 
of $\pi$ is a closed immersion, so the condition that the image is 
principal is indeed preserved under base change. On the other hand,
if we have $Z \in Z(\pi)$ such that $\pi(Z)\neq B$, but the image of
$B'$ is contained in $\pi(Z)$, then we observe that the first case of
condition (VIII) for $Z$ implies that $\sO_X(Y_Z) \otimes \sO_X(Y_Z^c)
\cong \sO_X$, which implies in turn that the second case of condition (VIII) 
is satisfied for $Z'$.

Additionally, it is
clear from the definition that the consistent condition is satisfied under
base change of a smoothing family. In particular, the
corollary implicitly includes the statement that in fibers
corresponding to smooth curves we recover the usual functor of linear
series.

Conceptually, we should think of independence of base change as being
a consequence of independence of limit linear series data. Indeed, in
simple cases (for instance, in a family with a unique singular fiber), 
there may be a natural minimal choice of a sufficient collection of 
concentrated multidegrees, but this will not be preserved under base
change: for instance, any collection of more than one concentrated 
multidegree will become redundant under restriction to a smooth fiber.
\end{remark}

Given the first statement of Corollary
\ref{cor:lls-independence}, we will henceforth drop the $\sD$ from
$\sG^r_d(X/B,\sD)$. We can formalize the functor $\sG^r_d(X/B)$ either
by fixing a choice of $\sD$ for each $X/B$, or by defining
$\sG^r_d(X/B)$ to be the limit over all choices of $\sD$ of the spaces
$\sG^r_d(X/B,\sD)$ under the isomorphisms provided by Corollary
\ref{cor:lls-independence}.  From either of these points of view, the
second statement of Corollary \ref{cor:lls-independence} then tells us
that if $X/B$ is a consistent smoothing family, and $B'\to B$ is any
morphism such that $X_{B'}\to B'$ is a (necessarily consistent) smoothing 
family, we have a canonical isomorphism
$$\beta_{B'/B}:\sG^r_d(X'/B')\risom\sG^r_d(X/B)\times_B B'.$$
Moreover, for any further morphism $B''\to B'$ such that
$X_{B''}\to B''$ is also a smoothing family, we have that
$$\beta_{B''/B}=\beta_{B'/B}|_{B''}\circ\beta_{B''/B'}.$$

\begin{notn}\label{notn:section-multideg}
Given a presmoothing family $\pi:X\to B$ with a
section $s:B\to X$, let $\Pic^{d,s}(X/B)$ denote the Picard scheme of
line bundles having degree $d$ on the component of each fiber of $\pi$
containing $s$, and degree $0$ on every other component.
\end{notn}

\begin{prop}\label{prop:lls-represent} If $\pi:X\to B$ is a
consistent smoothing family, then the functor $\sG^r_d(X/B)$ is
representable by a scheme $G^r_d(X/B)$ that is proper 
over $B$, and canonically compatible with base changes
that preserve the consistent smoothing family condition.
 
Moreover, we have that $G^r_d(X/B)$ has universal
relative dimension at least the Brill-Noether number 
$\rho$. If we further have that for some
$b \in B$, the space $G^r_d(X_{b})$ has dimension $\rho$, then we have
that $G^r_d(X/B)$ is Cohen-Macaulay and flat over $B$ at every point over 
$b$.
\end{prop}

The terminology of ``universal relative dimension at least
$\rho$'' is as introduced in Definition 3.1 of \cite{os21}.

\begin{proof} Given the definitions, representability is
rather standard. First, our functor is visibly a Zariski sheaf, so it
suffices to show representability locally on the base. Now,
we have already mentioned twisting by the
$\sO^Y$ allows us to move between any two multidegrees of fixed total
degree, and the hypothesized section means that line bundles of every
degree exist. Thus, if we choose any multidegree which is positive on 
every component of every fiber -- for instance, by summing all our
uniformly concentrated multidegrees -- we can produce a relatively ample 
line bundle. Passing to an open cover of $B$ if necessary we may assume we 
have a relatively ample divisor $D$.
The existence of the section also implies that the relative
Picard functor $\Pic^{d,s}(X/B)$ is representable, and carries a
Poincare line bundle $\widetilde{\sL}$. Twisting by a sufficiently
high multiple of $D$, we can then construct $G^r_d(X/B)$ as a closed
subscheme of a product of relative Grassmannians over
$\Pic^{d,s}(X/B)$, with the closed conditions given by vanishing along
the given multiple of $D$, intersected with the $(r+1)$st vanishing loci  
occurring in the definition of limit linear series.

Again passing to an open cover of the base, we
can assume our smoothing family is ``almost local'' in the sense of
Definition 2.2.2 of \cite{os20} (see also the proof of Proposition
\ref{prop:consistent-local}). We briefly sketch the ingredients of the
remaining argument before giving the details. In \cite{o-m1}, a theory
of `linked determinantal loci' is developed which is precisely tailored
to capture the conditions on $(r+1)$st vanishing loci for the maps
\eqref{eq:lls-map} which arise in the definition of limit linear series 
for a curve with two components. The main result of \cite{o-m1} is that 
these linked determinantal loci have the desired behavior in terms of 
codimension and Cohen-Macaulayness. In \cite{os25} and \cite{o-m1} it
is described how the linked determinantal locus construction can be
applied to deduce the desired results on limit linear series for arbitrary
curves of compact type (and in fact more generally). Thus, the desired
statements are essentially already contained in these papers. However,
there are some differences in hypotheses and details of constructions 
which should be addressed. The more basic is that the dimensionality
statement, proved as Theorem 6.1 of \cite{os25}, places more restrictions
on its families of curves, due to the desire to consider curves not of
compact type. However, the argument goes through verbatim in our case. 
The idea is that we can realize the limit linear series spaces for 
arbitrary curves of compact type as being cut out by conditions coming from 
pairs of adjacent components.

We then deduce Cohen-Macaulayness and 
flatness from Theorem 3.1 (see also Remark 3.5) of \cite{o-m1}. The
complication here is that there are several variants of our definition of limit
linear series used for different purposes, and while they all agree 
set-theoretically, we have not previously shown that their scheme structures
agree. Specifically, the scheme structure used in \cite{o-m1} imposes the 
condition on the $(r+1)$st vanishing
locus of \eqref{eq:lls-map} on a smaller collection of
multidegrees; thus, what we know is that our scheme structure is a closed 
subscheme of the one in \cite{o-m1}, with the same support, and in order
to complete the proof of the proposition we want to
show that it is in fact the same scheme. While we consider all multidegrees
$\md$ in our definition, the definition used in \cite{o-m1} considers only
multidegrees which are nonnegative on all components and equal to zero on
all but two adjacent components. We will show that the schemes agree on
fibers; it will then follow from Nakayama's lemma and the flatness of the 
larger scheme that the two scheme structures must agree. 
Now, the argument of the penultimate paragraph of the proof of Theorem 4.3.4 
of \cite{os20} shows that imposing the vanishing conditions on the 
multidegrees considered in \cite{o-m1} then implies that the same 
conditions are satisfied on all multidegrees which are ``nonnegative''
in the sense of being nonnegative on every $Y \in Y(\pi)$. Then, the 
proof of Proposition 3.4.12 of \cite{os20} 
shows that for any multidegree $\md$, there is some such nonnegative 
multidegree $\md'$ such that the kernel of \eqref{eq:lls-map} for $\md'$
injects into the kernel for $\md$; it then follows from Corollary
B.3.5 of \cite{os20} that the vanishing condition for $\md'$ implies
the same vanishing condition for $\md$. 
We thus obtain the desired agreement of scheme structures.
\end{proof}

\begin{cor}\label{cor:lls-descend} Suppose $\pi:X \to B$ is a
presmoothing family, and we also fix $r,d$. Then there exists a proper
algebraic space $G^r_d(X/B)$ over $B$ such that for any base change
$B' \to B$ making the resulting $\pi':X' \to B'$ into a consistent
smoothing family, we have a functorial identification $G^r_d(X/B)
\times _B B' \cong G^r_d(X'/B')$.

If moreover $\pi$ admits a section $s$, then
$G^r_d(X/B)$ maps to $\Pic^{d,s}(X/B)$, compatibly with the above
identification.

In general, $G^r_d(X/B)$ has universal relative
dimension at least the Brill-Noether number $\rho$, and if for some 
$b$ we have that the fiber
of $G^r_d(X/B)$ has dimension $\rho$ over $b$, then $G^r_d(X/B)$
is Cohen-Macaulay and flat over $b$.
\end{cor}

\begin{proof}
  Let $T$ be the full subcategory of
  the big \'etale site of $\Sch_B$ consisting of arrows $B''\to B$ such that
  $X_{B''}/B''$ is a consistent
  smoothing family. By Propositions \ref{prop:presmooth-smooth} and
  \ref{prop:consistent-local}, 
  there is an \'etale cover $B'\to B$ such that $X':=X\times_B B'$ is a 
  consistent smoothing family. It follows from Remark
  \ref{rem:descent-lemma-2} that $T$ satisfies the 
  conditions of Lemma \ref{lem:descent-daddy}. Applying that
  Lemma, there is a sheaf
$\sG^r_d(X/B)$ on the big \'etale site of $\Sch_B$ whose value on any
$B''\to B$ over which $X_{B''}/B''$ is a consistent smoothing family is the
sheaf $\sG^r_d(X_{B''}/B'')$ of Definition \ref{defn:d-lls}. Since
$\sG^r_d(X_{B'}/B')$ is a proper scheme,
we have that $\sG^r_d(X/B)$ is itself a
proper algebraic space over $B$ (see Tags 083R, 0410,
03KG and 03KM of \cite{stacks-proj}).  

Suppose $\pi$ admits a section $s$. For any $B'' \to B$ such that $X_{B''}/B''$
is a consistent smoothing family, there is a forgetful morphism
$\sG^r_d(X_{B''}/B'')\to\Pic^{d,s}(X_{B''}/B'')$ defined by choosing the 
(unique) twist 
of the underlying line bundle that has multidegree concentrated along the
image of $s$, as described in Notation
\ref{notn:section-multideg}. This defines a morphism of sheaves on
$T$. By Lemma \ref{lem:descent-daddy}, these extend uniquely to a
morphism $\sG^r_d(X/B)\to\Pic^{d,s}(X/B)$ of sheaves on $\Sch_B$.
\margh{standardize base change/restriction notn/terms}

Universal relative dimension for stacks is defined
(Definition 7.1 of \cite{os21}) in terms of descent from a smooth
cover, so the desired statement follows from Proposition
\ref{prop:lls-represent}. Similarly, since fiber dimension is
invariant under base extension and flatness and Cohen-Macaulayness descend, 
the final statements also follow from Proposition \ref{prop:lls-represent}.
\end{proof}


More generally, we have the following.

\begin{cor}\label{cor:lls-descend-stack} Suppose $\cB$ is a
Deligne-Mumford stack and $\widetilde{\pi}:\cX \to \cB$ is a curve,
and there exists an \'etale cover $B \to \cB$ by a scheme such that the
induced $\pi':X \to B$ is a presmoothing family. Then for any fixed
$r,d$, there exists a Deligne-Mumford stack $\cG^r_d(\cX/\cB)$ which
is proper and a relative algebraic space over $\cB$ such that for any
$B' \to \cB$ making the resulting $\pi':X' \to B'$ into a consistent
smoothing family, we have a functorial identification
$\cG^r_d(\cX/\cB) \times _{\cB} B' \cong G^r_d(X'/B')$.

If moreover $\pi$ admits a section $s$, then
$\cG^r_d(\cX/\cB)$ maps to $\Pic^{d,s}(\cX/\cB)$, compatibly with the
above identification.

In general, $\cG^r_d(X/B)$ has universal relative dimension at least 
the Brill-Noether number $\rho$, and if for some $b$ we have that the fiber
of $\cG^r_d(X/B)$ has dimension $\rho$ over $b$, then $\cG^r_d(X/B)$
is Cohen-Macaulay and flat over $b$.
\end{cor}

\begin{proof} The proof is identical to the proof of Corollary
\ref{cor:lls-descend}, using the big \'{e}tale site of $\cB$ in place
of $\Sch_B$.
\end{proof}

In view of Proposition \ref{prop:presmooth-smooth}, the following is an 
immediate consequence of Corollary \ref{cor:lls-descend-stack}.   

\begin{cor}\label{cor:univy} There is a proper relative
algebraic space
$$\cG^r_d\to\overline{\cM}_{g}^{\ct}$$
with the property that for any morphism $B\to\overline{\cM}_{g}^{\ct}$ such
that the pullback of the universal family to $B$ is a consistent
smoothing family $X\to B$, there is a natural isomorphism
$$\cG^r_d\times_{\overline{\cM}_{g}^{\ct}} B\cong G^r_d(X/B)$$
of $B$-spaces.
\end{cor}

We have thus constructed universal moduli stacks of limit
linear series, proving Theorem \ref{thm:universal}. 

Finally, we want a good description of the fibers of our
construction.  In particular, we have the following.

\begin{cor}\label{cor:lls-fibers} Suppose that $k$ is a field,
$X_0$ a curve of compact type over $k$, and $k'/k$ is a finite
Galois extension over which $X_0$ is totally split.  
Given $r,d$, we have the scheme $G^r_d(X'_0/k')$ and the algebraic
space $G^r_d(X_0/k)$, and the $k$-points of $G^r_d(X_0/k)$ are
canonically identified with the $k'$-points of $G^r_d(X'_0/k')$ which
are invariant under the natural $\Gal(k'/k)$-action.
\end{cor}
\begin{proof}
  This follows from Remark \ref{rem:descent-lemma-3}, since the
  category $T_0$ (as described in Corollary \ref{cor:lls-descend})
  can be chosen to consist of the iterated fiber products of $\Spec
  k'$ over $\Spec k$.
\end{proof}

\section{Applications over real and $p$-adic fields}\label{sec:real}

We now apply our results to study rational points over real
and $p$-adic fields, in particular proving Theorem \ref{thm:lls-real}
and Corollary \ref{cor:lls-real} in the real case. We first prove the
following basic result. (While we believe this should be in the
literature somewhere with an actual proof, we have failed to find it.)

\begin{prop}\label{prop:ratl-constant} Suppose given a field
$K$ which is either $\RR$ or a $p$-adic field, a smooth curve $B$ over
$K$, a point $b \in B(K)$, and a finite flat morphism $\pi:Y \to B$ of
$K$-schemes. Suppose further that no ramification point of $\pi$ over
$b$ is a $K$-point.  Then there is an analytic neighborhood $U$ of $b$
such that for all $b' \in U(K)$, none of the $K$-rational points in
the fibers $Y_{b'}$ are ramified, and the number of $K$-rational points
over $U$ is constant.
\end{prop}

In the above, $B(K)$ and $Y(K)$ denote sets of $K$-rational
points.

\begin{proof} Let $\nu:(U,u)\to (B,b)$ be an \'etale
neighborhood of $b$, with $u \in U(K)$, over which $Y$ splits as a disjoint 
union
$Y\times_B U=\bigsqcup_{i=1}^n Y_i$, with $Y_i\to U$ an isomorphism
for $i=1,\ldots,m$ and $Y_j\to U$ a morphism whose fiber over $b$
contains no $K$-points for $j=m+1,\dots,n$. Since $\nu$ is an 
analytic-local isomorphism,
it suffices to prove the result for $Y_U\to U$ and thus we may assume
that $Y$ breaks up in the manner described. We will thus replace $B$
by $U$ and work under the additional assumption that $Y$ decomposes.

The result follows from the following key property: 
the induced map $Y_i(K)\to B(K)$ is a proper map of Hausdorff
topological spaces. To prove this it suffices (by choosing embeddings and
using the fact that the analytic topology is finer than the Zariski
topology) to prove the same statement under the assumption that
$B=\AA^n$ and $Y_i=\PP^m\times\AA^n$, the map $Y_i\to B$ being replaced
by the second projection map.

We wish to show that the
preimage of a compact set $S$ in $\AA^n(K)$ is compact. But this preimage
is just $S\times\PP^m(K)$. Since a product of two compact spaces is
compact, it thus suffices to show that $\PP^n(K)$ is itself
compact. To prove this, it suffices to show that $\PP^n(K)$ is the
image of a compact set under a continuous map. If $K=\RR$, we know that
$\PP^n(K)$ is the image of the unit sphere in $\AA^{n+1}(K)$. If $K$
is non-archimedean, write $\cO$ for the ring of integers.
Let $C_i\subset\cO^{n+1}$ be the subset for which the $i$th
coordinate is $1$. Since $\cO$ is compact, $C_i$ is itself
compact. Scaling a point of $\PP^n(K)$ by the inverse of the coordinate with the
largest absolute value, we see that there is a surjection $$\bigsqcup_{i=1}^{n}
C_i\to\PP^n(K).$$ This establishes the assertion.

With this property in hand, we can conclude the proof.
Since the preimage $y_j$ of $b$
in $Y_j$ for $j=m+1,\ldots,n$ is not a $K$-point, each $Y_j$ for $j$
in this range has the property that $\pi(Y_j(K))\subset B(K)$ does not
contain $b$, so that $B(K)\setminus\pi(Y_j(K))$ is an open subset not
containing $b$, which means we can shrink $B$ so that $Y_j(K)$ is
empty for each $j=m+1,\ldots,n$. Since each other map $Y_j\to B$,
$j=1,\ldots,m$ is an isomorphism, we get the desired result.
\end{proof}

We next recall the Eisenbud-Harris definitions of limit linear
series.

\begin{defn}\label{defn:eh-lls} Let $X_0$ be a totally split curve of
compact type, with dual graph $\Gamma$. For $v \in V(\Gamma)$, let
$Y_v$ be the corresponding component of $X_0$, and for $e \in
E(\Gamma)$, let $Z_e$ be the corresponding node. An
\emph{Eisenbud-Harris limit} $\fg^r_d$ on $X_0$ is a tuple
$(\sL^v,V^v)_{v \in V(\Gamma)}$ of $\fg^r_d$s on the $Y_v$, satisfying
the condition that for each node $Z_e$, if $Y_v,Y_{v'}$ are the
components containing $Z_e$, then for $j=0,\dots,r$, we have
\begin{equation}\label{eq:eh-ineq} a_j+a'_{r-j} \geq
d,
\end{equation} where $a_0,\dots,a_r$ and
$a'_0,\dots,a'_r$ are the vanishing sequences at $Z_e$ of
$(\sL^v,V^v)$ and $(\sL^{v'},V^{v'})$ respectively.

The limit $\fg^r_d$ given by $(\sL^v,V^v)_v$ is said
to be \emph{refined} if \eqref{eq:eh-ineq} is an equality for all $e$
and $j$.
\end{defn}

By definition, the set of Eisenbud-Harris limit linear series
is a subset of the product of the spaces $G^r_d(Y_v)$ of usual linear
series on the $Y_v$.  If we choose for each $e \in E(\Gamma)$
sequences $a_{\bullet},a'_{\bullet}$ for which \eqref{eq:eh-ineq} is
satisfied with equality, we obtain a product of spaces of linear
series on the smooth curves $Y_v$ with imposed ramification at the
$Z_e$, and taking the union over varying choices of the
$a_{\bullet},a'_{\bullet}$, one sees that we obtain every
Eisenbud-Harris limit linear series on $X_0$. This then induces a
natural scheme structure:

\begin{notn} In the situation of Definition \ref{defn:eh-lls},
let $G^{r,\EH}_d(X_0)$ denote the moduli scheme of Eisenbud-Harris
limit linear series, with scheme structure as described above.
\end{notn}

In order to state the strongest comparison results, the following
preliminary proposition will be helpful.

\begin{prop}\label{prop:refined-dense} Let $X_0$ be a totally split curve of
compact type, and given $r,d$, suppose that the moduli space
$G^{r,\EH}_d(X_0)$ of limit $\fg^r_d$s has dimension equal to the
Brill-Noether number $\rho$. Then the
refined limit linear series are dense in $G^{r,\EH}_d(X_0)$.
\end{prop}

\begin{proof} We will prove the following more precise statement:
Suppose that $G^{r,\EH}_d(X_0)$ has dimension $\rho$, and
we are given, for each pair $(e,v)$ of an edge and adjacent vertex in
$\Gamma$, a sequence 
$$0 \leq a^{(e,v)}_0 < a^{(e,v)}_1 < \dots < a^{(e,v)}_r \leq d,$$
such that for every edge $e$, connecting $v$ to $v'$, we have
\begin{equation}\label{eq:eh-ineq-2}
a^{(e,v)}_j+a^{(e,v')}_{r-j} \geq d\text{ for } j=0,\dots,r.
\end{equation}
Then the closed subset $S_{(a^{(e,v)})}$ of $G^{r,\EH}_d(X_0)$ consisting 
of limit linear
series $(\sL^v,V^v)_v$ with $V^v$ having vanishing sequence at $Z_e$ 
at least equal to $a^{(e,v)}$ for all $(e,v)$ has dimension
\begin{equation}\label{eq:rho-codim}
\rho-\sum_{e,j} \left(a^{(e,v)}_j+a^{(e,v')}_{r-j} - d\right)
\end{equation}
if it is nonempty. From this, it will follow immediately that the refined 
limit linear series are dense.

Arguing the contrapositive, suppose there exist sequences $a^{(e,v)}$
as above such that $S_{(a^{(e,v)})}$ has dimension strictly bigger than
asserted; we claim that if \eqref{eq:eh-ineq-2} is strict for any
$(e,v)$ and $j$, then we can decrease some $a^{(e,v)}_j$ while
preserving the condition that the resulting 
$S_{(a^{(e,v)})}$ has too large a dimension. Now, the subset
$S_{(a^{(e,v)})}$ is by definition a product of spaces of linear series
with imposed ramification on the $Y_v$; if we denote by $\rho_v$ and
$d_v$ the expected and actual dimensions of these spaces respectively,
we always have $d_v \geq \rho_v$ for all $v$, and $\sum_v \rho_v$ is
equal to the expression in \eqref{eq:rho-codim}.
We therefore have that
the condition that $S_{(a^{(e,v)})}$ has too large dimension is equivalent
to saying that $d_v>\rho_v$ for some $v$. Now, suppose that we have
$(e,v)$ and $j$ with strict inequality in \eqref{eq:eh-ineq-2}, and
let $v'$ be the other vertex adjacent to $e$. If $d_{v''}>\rho_{v''}$
for some $v'' \neq v$, then we can replace $a^{(e,v)}_j$ with
$d-a^{(e,v')}_{r-j}$ for $j=0,\dots,r$, which only changes $d_v$ and
$\rho_{v''}$, and in particular preserves that $d_{v''}>\rho_{v''}$.
On the other hand, if $d_{v''}=\rho_{v''}$ for all $v'' \neq v$, then
necessarily $d_v>\rho_v$, and we can similarly replace 
$a^{(e,v')}_{r-j}$ with $d-a^{(e,v)}_{j}$ for $j=0,\dots,r$, preserving
the hypothesis that $d_v>\rho_v$. This proves the claim, and
iterating this procedure, we eventually produce a subset of 
$G^{r,\EH}_d(X_0)$ of dimension greater than $\rho$. This proves the
desired contrapositive statement.
\end{proof}

The most basic comparison result is then the following.

\begin{lem}\label{lem:lls-agree-sets} In the situation of
Definition \ref{defn:eh-lls}, choose the sufficient collection of
concentrated multidegrees with $I=V(\Gamma)$, and $\md_v$ having
degree $d$ on $Y_v$ and degree $0$ on the other components. Then
restriction from $X_0$ to $Y_v$ induces a bijection from $G^r_d(X_0)$
to $G^{r,\EH}_d(X_0)$.
\end{lem}

Indeed, this is the rank-$1$ case of Lemma 4.1.6 and Proposition 4.2.9 of
\cite{os20}.

In general, it is not clear that the scheme structures on
$G^r_d(X_0)$ and $G^{r,\EH}_d(X_0)$ agree, but it follows from
\cite{os20} and \cite{o-m1} that under the most common circumstances
we will have agreement. The improved statement using Proposition
\ref{prop:refined-dense} is as follows.

\begin{thm}\label{thm:lls-agree} The map $G^r_d(X_0) \to
G^{r,\EH}_d(X_0)$ constructed in Lemma \ref{lem:lls-agree-sets} is an
isomorphism of schemes when the following conditions are satisfied:
\begin{enumerate} 
\item[(I)] $G^{r,\EH}_d(X_0)$ has dimension equal to the Brill-Noether number
$\rho$; 
\item[(II)] either $\rho=0$, or $G^{r,\EH}_d(X_0)$ is reduced.
\end{enumerate}
\end{thm}

\begin{proof} Proposition \ref{prop:refined-dense} implies that under
condition (I), we also have that the refined limit linear series are
dense. The case that $G^{r,\EH}_d(X_0)$ is reduced is then
Corollary 3.3 of \cite{o-m1}. In the case $\rho=0$, density implies 
that in fact every limit linear series is refined, in which case
the desired statement is Proposition 4.2.6 of \cite{os20}.
\end{proof}

We can now give the proof of our main theorem on limit linear
series over local fields.

\begin{proof}[Proof of Theorem \ref{thm:lls-real}] First, by
Corollary \ref{cor:lls-descend}, we may restrict to a neighborhood of
$b_0$ over which $G^r_d(X/B)$ is quasi-finite and flat; we then have that 
it is finite (see Tag 0A4X of \cite{stacks-proj}), and then it immediately
follows that it is also a scheme.  By Theorem \ref{thm:lls-agree},
our hypotheses imply that $G^r_d((X_0)_L) \cong G^{r,\EH}_d((X_0)_L)$,
and we can then invoke Corollary \ref{cor:lls-fibers} to conclude that
$G^r_d(X_0)$ is finite, that $G^r_d(X_0)(K)$ has size exactly $n$, and that the 
the scheme structure at each element of $G^r_d(X_0)(K)$ is reduced.  The 
desired statement then follows from
Proposition \ref{prop:ratl-constant}.
\end{proof}

Following the argument of \cite{os2} and \cite{os4}, we can
now reduce the question of existence of real curves with given numbers
of real linear series to statements on real points of intersections of
Schubert cycles.  Specifically, we have the following.

\begin{cor}\label{cor:lls-schubert} Fix $r<d$, and suppose $K$
is either $\RR$ or a $p$-adic field. Given a complete flag
$$(0)=F_0 \subseteq F_1 \subseteq \dots \subseteq F_{d+1}=K^{d+1},$$ 
let $\Sigma(F_{\bullet})$ be the Schubert cycle in
$G(r+1,d+1)$ consisting of $(r+1)$-dimensional subspaces meeting
$F_{d-1}$ in codimension $1$ (inside the subspace). 
Given $P \in \PP^1$, let $F^P_{\bullet}$
denote the flag in $\Gamma(\PP^1,\sO(d))$ determined by order of vanishing
at $P$. Given $g$ with $g=(r+1)(r+g-d)$, suppose that there exist
$P_1,\dots,P_g \in \PP^1(\bar{K})$ such that the divisor
$P_1+\dots+P_g$ is $\Gal(\bar{K}/K)$-invariant, and with
$$\bigcap_{i=1}^g \Sigma(F^{P_i}_{\bullet})$$
finite, having $n$ $K$-rational points, all of which
are reduced.

Then there exist smooth projective curves $X$ of genus
$g$ over $K$ such $G^r_d(X)$ is finite, with exactly $n$ $K$-rational
points, all of which are reduced.
\end{cor}

\begin{proof} Making use of our new machinery, the proof is closely
based on the proof of Theorem 2.5 of \cite{os2}.
Given the Galois invariance of $P_1+\dots+P_g$,
we can construct a (not necessarily totally split) curve $X_0$ over $K$ of
compact type by attaching suitable elliptic tails to the
$P_i$. Letting $L/K$ be a finite Galois extension over which all the
$P_i$ are rational, we will have that $(X_0)_L$ is a totally split
curve of compact type. Now, as described in the proof of Theorem 2.5
of \cite{os2}, the only possible vanishing sequences at the nodes in
the $\rho=0$ case are $d-r-1,d-r,\dots,d-2,d$ on the elliptic tails,
with complementary sequence $0,2,3,\dots,r+1$ on the rational ``main
component.'' Moreover, there is a unique
linear series on each elliptic tail with the desired vanishing at
$P_i$, even scheme-theoretically, so we find that
$G^{r,\EH}_d((X_0)_L)$ consists entirely of refined limit linear
series, and is isomorphic to the space of $\fg^r_d$s on $\PP^1$ with
vanishing sequence $0,2,3,\dots,r+1$ at each $P_i$. This latter space
is precisely $\bigcap_{i=1}^g \Sigma(F^{P_i}_{\bullet})$ (considered over
$L$), since $F^{P_i}_{d-1}$ is the space of polynomials vanishing to order
at least $2$ at $P_i$. By construction,
$\Gal(L/K)$ acts on $G^{r,\EH}_d((X_0)_L)$ via its action on $X_0$
itself, which means that our isomorphism is $\Gal(L/K)$-equivariant,
and the invariant points on both sides are identified. Thus, the
hypotheses of Theorem \ref{thm:lls-real} are satisfied, and taking any
presmoothing family with special fiber $X_0$ (see \S 2 of \cite{os4} for
justification that this is possible), by considering smooth fibers $X$
sufficiently near $X_0$, we obtain the desired properties for
$G^r_d(X)$.
\end{proof}

We can now apply Corollary \ref{cor:lls-schubert} to examples
of Eremenko and Gabrielov in the real case.

\begin{proof}[Proof of Corollary \ref{cor:lls-real}] Using
Corollary \ref{cor:lls-schubert}, this is immediate from Examples 1.3
and 2.5 of \cite{e-g3}, where their $m$ is our $d-1$.  
\end{proof}

\begin{rem} In terms of producing real curves with few real $\fg^1_d$s,
the result of Corollary \ref{cor:lls-real} is the best possible using
our given degeneration, since Eremenko and Gabrielov prove that their
examples yield the smallest possible number of real points in the 
corresponding intersection of Schubert cycles. Nonetheless, it is 
\emph{a priori} possible (in the case $d$ is even) that there exist real 
curves with fewer real $\fg^1_d$s, but that these are not ``close enough''
to the degenerate curves we consider to be studied via our techniques.
\end{rem}

\begin{rem} It is an interesting phenomenon that the examples
of Eremenko and Gabrielov achieving the minimum possible number of
real solutions (in the Schubert cycle setting) occur not with all 
ramification points non-real, but
with exactly two real. Computations of Hein, Hillar and Sottile 
\cite{h-h-s1} in small degree
appear to show that on the other hand, the numbers obtained by Cools
and Coppens are the minimum achievable when all the ramification
points are non-real.  

In fact, every such explicit Schubert calculus example also gives an
example where we can produce higher-genus curves having the same number
of real $\fg^1_d$s, using Corollary \ref{cor:lls-schubert}.
However, we are not aware of infinite families of examples other than
those of \cite{e-g3}.
\end{rem}

\bibliographystyle{amsalpha} 
\newcommand{\noopsort}[1]{} \newcommand{\printfirst}[2]{#1}
  \newcommand{\singleletter}[1]{#1} \newcommand{\switchargs}[2]{#2#1}
\providecommand{\bysame}{\leavevmode\hbox to3em{\hrulefill}\thinspace}
\providecommand{\MR}{\relax\ifhmode\unskip\space\fi MR }
\providecommand{\MRhref}[2]{%
  \href{http://www.ams.org/mathscinet-getitem?mr=#1}{#2}
}
\providecommand{\href}[2]{#2}

\end{document}